# Unbiased Simulation of Rare Events in Continuous Time

**James Hodgson¹ · Adam M. Johansen¹,³ · Murray Pollock²,³**

**Abstract**

For rare events described in terms of Markov processes, truly unbiased estimation of the rare event probability generally requires the avoidance of numerical approximations of the Markov process. Recent work in the exact and $\varepsilon$-strong simulation of diffusions, which can be used to almost surely constrain sample paths to a given tolerance, suggests one way to do this. We specify how such algorithms can be combined with the classical multilevel splitting method for rare event simulation. This provides unbiased estimations of the probability in question. We discuss the practical feasibility of the algorithm with reference to existing $\varepsilon$-strong methods and provide proof-of-concept numerical examples.

## 1 Introduction

*Rare* events are those which have (very) low probability of occurrence. Estimating the probability of rare events is important, among other places, throughout the natural and social sciences; see, for example, [Part II] Rubino and Tuffin (2009) for a broad range of applications. The case of interest here is that where the rare event corresponds to a continuous-time Markov process hitting a particular set before it enters some other recurrent

✉ James Hodgson
james.hodgson@warwick.ac.uk

Adam M. Johansen
a.m.johansen@warwick.ac.uk

Murray Pollock
murray.pollock@newcastle.ac.uk

¹ Department of Statistics, University of Warwick, CV4 7AL Coventry, UK

² School of Mathematics, Statistics and Physics, Newcastle University,
NE1 7RU Newcastle-upon-Tyne, UK

³ The Alan Turing Institute, The British Library, London NW1 2DB, UK

set. This setting has attracted considerable attention in the literature, and general solutions centre around simulation-based methods.

The principal approaches fall into two broad categories, *importance sampling* and *splitting*. In importance sampling, one simulates from a process for which the event of interest is more likely to occur, and corrects for the change of sampling distribution using importance weights. With splitting methods, trajectories which approach the rare set (in an appropriate sense) are replicated to allow lower-variance estimation of the target probability. This paper is concerned with splitting methods, in particular with implementing such methods with no bias for a broad class of continuous-time processes. Existing methods depend upon time-discretisation and hence introduce a difficult to quantify bias. It is shown in this article that the adaptation of ideas from the field of $\epsilon$-strong simulation to this context allows this bias to be avoided. Let $(X(t) : t \geq 0)$ be a continuous-time Markov process in $\mathbb{R}^d$, and let $A, B \subset \mathbb{R}^d$ be disjoint sets, with $A$ positive recurrent for $X$. The problem of interest is that of efficiently estimating the probability that $X$ reaches set $B$ before set $A$ when this probability is very small. That is, writing $\tau_S$ for the first hitting time of a set $S$, the objective is to estimate

$$p = \mathbb{P}(\tau_B < \tau_A) \ll 1.$$

The assumption that $p \ll 1$ rules out direct Monte Carlo estimation, since the computational cost of generating the event $\{\tau_B < \tau_A\}$ enough times to get a reliable estimate will be impractically high. In fact, the relative variance of the naive estimator obtained from $N$ Monte Carlo simulations is $p(1-p)/N \cdot p^{-2} \approx 1/(Np)$, which suggests that a large multiple of $p^{-1} \gg 1$ trials is needed for a reasonable estimate.

*Multilevel splitting* (MLS) is a popular algorithm based on targeting the rare event via a sequence of more likely events. The idea goes back to the 1951 paper of Kahn and Harris Kahn and Harris (1951) (who in turn attribute the idea to von Neumann), which discusses an application to the transmission of particles through an impeding barrier in the context of nuclear shielding. The method is to choose a sequence of nested sets $B_1 \supset B_2 \supset \cdots \supset B_m = B$, all disjoint from $A$, and use a particle system to sequentially estimate $\mathbb{P}(\tau_{B_i} < \tau_A)$. Starting with a particle system of a large enough size, $N$, a reasonable fraction will reach $B_1$, allowing an estimate of $p_1$. Then, by branching (or "splitting") those which do into $R_i$ copies, a healthy population can be maintained to estimate the subsequent probabilities. A detailed presentation is given in Sect. 2.1.

Splitting algorithms have been independently rediscovered many times and in many variants since the work of Kahn and Harris (1951). Prominent examples include the *repetitive simulation trials after reaching thresholds* (RESTART) algorithm of Villén-Altamirano and Villén-Altamirano (1994), developed for modelling packet loss probabilities in telecommunications, and the *pruning-enriched Rosenbluth method* (PERM) of Grassberger (1997) for simulating polymer chains.

The MLS algorithm we present in Sect. 2.1 is that found in Garvels (2000) which addresses various implementation issues such as the choice of levels and importance function. The unbiasedness of the algorithm for discrete-time processes is shown rigorously in Amrein and Künsch (2011), which identifies and resolves an issue in the original argument of Garvels (2000). The construction of confidence intervals and the optimal choice of tuning parameters under cost constraints are addressed in Lagnoux-Renaudie (2006, 2008). A characterization of the asymptotic properties of this algorithm, including a central limit theorem, are given in Del Moral and Lezaud (2006).

Choosing the nested sets and other parameters of the algorithm to maintain a particle population of stable size, rather than one which dies out or explodes, can be difficult. One practical variant which removes the difficulty of choosing the *splitting ratios* $R_i$ in advance is that of Lagnoux-Renaudie (2009), in which a first particle system is used to estimate the $R_i$, and a second system uses these estimated values to estimate $p$. An alternative idea is to construct the levels $B_i$ adaptively, for example via the scheme of Cérou and Guyader (2007). A generalisation of this scheme has recently been shown be unbiased in Bréhier et al. (2016).

The proof of unbiasedness in Amrein and Künsch (2011) also holds for a variant in which the initial system of $N$ particles is kept at fixed size by sampling new trajectories uniformly at random (with replacement) from the surviving trajectories at each level. This variant is also discussed in Garvels (2000) under the name of fixed splitting. It is a type of Sequential Monte Carlo method and can be understood within the framework of Del Moral (2004). In this version the sets $B_i$ are still chosen in advance, but the number of particles is fixed at $N$ for the duration of the algorithm. Rather than independently "splitting" each path which survives to $B_i$ into a pre-determined number of offspring, exactly $N$ particles are resampled (i.e. sampled with replacement) from among the surviving particles. This removes the difficulty of choosing a suitable splitting ratio in order to arrive at a stable population size, but it is more difficult to understand the variance properties of this algorithm and even the asymptotic variance expression is somewhat more complex than the variance of the simple algorithm.

As discussed, in principle at least the estimate obtained using MLS is unbiased. But in order to implement the algorithm, it is necessary to simulate hitting times and locations of the sets $B_i$ for the process $X$. When dealing with discrete-time processes, there is often no difficulty in representing a full sample path, and so a direct implementation gives unbiased estimates. However, when the process of interest evolves in continuous-time, Monte Carlo sample paths are usually constructed using a discrete-time numerical scheme, causing bias in the resulting estimates—at least outside the setting of finite activity pure jump processes. We discuss this issue and some related literature in Sect. 2.2. This paper introduces a framework which combines algorithms for the *$\varepsilon$-strong simulation* of sample paths of $X$ (sometimes termed tolerance-enforced simulation) with MLS in order to obtain a truly unbiased algorithm.

For a chosen time horizon $[s, t]$, $\varepsilon$-strong methods are explicit constructions of a family of random processes $\tilde{X}^{\varepsilon}[s, t]$ indexed by the parameter $\varepsilon > 0$. These processes are defined jointly on the same probability space as $X$, admit a finite-dimensional representation, and satisfy

$$\sup_{r \in [s,t]} \|X(r) - \tilde{X}^{\varepsilon}(r)\| < \varepsilon$$

almost surely, for an appropriate norm $\|\cdot\|$, typically the supremum norm $\|\cdot\|_{\infty}$. Such paths constrain, almost surely, the range of $X$ over the time interval $[s, t]$ to within the tolerance $\varepsilon$. Moreover, for $\eta < \varepsilon$, such schemes allow one to sample $\tilde{X}^{\eta}$ conditional on $\tilde{X}^{\varepsilon}$. This means it is possible to sample, exactly, Bernoulli random variables which indicate whether a path of $X$ entering given subsets of $\mathbb{R}^d$. We use such $\varepsilon$-strong samplers to construct a modified procedure for MLS, and establish that this modified procedure is unbiased — and, in contrast to the abstract MLS algorithm without time discretisation, can be implemented for a class of continuous-time processes.

The literature on $\varepsilon$-strong simulation is closely related to that on the *exact simulation of diffusions*, (Beskos and Roberts, 2005; Beskos et al, 2006). Both constructions are rejection samplers on the space of continuous paths, using Brownian-like proposals, in which the acceptance probability for a given path can be assessed by looking at only finitely many points. The incorporation of exact simulation into an SMC algorithm has been successfully carried out in the filtering setting (Fearnhead et al, 2008); our approach here is different, using $\varepsilon$-strong simulation algorithms within the SMC context to allow for unbiased estimation of rare event probabilities. The unbiasedness of our approach follows from the possibility of exactly determining whether sample paths cross particular barriers, rather from the use of unbiased random weights as in the *random weight particle filter* setting of Fearnhead et al. (2008).

The first $\varepsilon$-strong algorithm, for Brownian motion, was given in Beskos et al. (2012). The construction exploits a certain representation of the escape probability of one-dimensional Brownian motion from a bounded interval. The contributions of Pollock et al. (2016) showed that exact simulation was possible for a much wider class of diffusions and jump-diffusions than merely Brownian motion (which include multidimensional diffusions amenable to Lamperti transformation, see Lamperti (1964)). This precise construction is what we use for the development of our approach applied to MLS.

Although we focus in our implementation upon the algorithmic construction of Pollock et al. (2016), our framework can employ any $\varepsilon$-strong approach with certain basic properties and we note that this area is being actively researched and a number of other such schemes can be found in the literature. Thus, in principle, the approach developed in this paper provides a mechanism for the unbiased estimation of rare event probabilities associated with a broad and increasing category of stochastic processes. A more general multidimensional $\varepsilon$-strong construction appears in Blanchet et al. (2017), which employs altogether different techniques inspired by the theory of rough paths. Rather than sampling the diffusion path itself exactly, one instead discretely samples the Brownian path driving the diffusion. A rough path-type continuity result guarantees that passing this discrete Brownian sample through a modified Euler scheme gives a sample within $\varepsilon$-tolerance of the diffusion path. Chen et al. (2019) contains a related construction for fractional Brownian motion of any Hurst index, and for solutions to SDEs driven by fractional Brownian motion with Hurst index $H > \frac{1}{2}$. An $\varepsilon$-strong algorithm for a process which is not the solution to an SDE is described in Cázares et al. (2019). Recent applications of this methodology include Mider et al. (2019).

This paper will demonstrate that the information provided by $\varepsilon$-strong algorithms can be exploited within MLS algorithms to obtain an exact estimate of the desired rare event probability. In Sect. 2, we introduce multilevel splitting formally, discuss the problem of estimation bias into practical implementations, and present $\varepsilon$-strong sampling as a way to overcome this. In Sect. 3, we develop in detail a method for implementing epsilon-strong samplers inside MLS algorithms and establish the unbiasedness of the resulting estimates. In Sect. 4, we offer some simple numerical simulations to illustrate the method. We conclude with a brief discussion in Sect. 5.

## 2 Background

This section describes the existing work upon which the developments of Sect. 3 depend. In Sect. 2.1, we give a full description of multilevel splitting and its SMC variant, together with discussion of the unbiasedness of the corresponding estimators. In Sect. 2.2, we explain why unbiasedness generally does not hold in any existing implementation of these methods. Finally, Sect. 2.3 introduces a formalism of $\varepsilon$-strong simulation (following Blanchet et al. (2017)), and describes how it can be used to address the shortcomings detailed in Sect. 2.2.

## 2.1 Multilevel Splitting

Recall that multilevel splitting requires the specification of a sequence of nested events $\mathbb{R}^d \supset B_1 \supset B_2 \supset \cdots \supset B_m = B$ in such a way that the probabilities

$$p_i := \begin{cases} \mathbb{P}(\tau_{B_1} < \tau_A), & i = 1 \\ \mathbb{P}(\tau_{B_i} < \tau_A \mid \tau_{B_{i-1}} < \tau_A), & 2 \leq i \leq m \end{cases}$$

are large relative to $p$, and may consequently be estimated more efficiently.

In order to do this, it is convenient to assume the existence of a continuous function $\xi : \mathbb{R}^d \to \mathbb{R}$ of which the boundaries of $A$ and $B_i$ are *level sets*. That is, we suppose that there are real numbers $z_A < z_1 < z_2 < \cdots < z_m = z_B$ such that

$$A = \xi^{-1}\left((-\infty, z_A]\right), \ B_i = \xi^{-1}\left([z_i, \infty)\right),$$

with boundaries $\partial A = \xi^{-1}(\{z_A\}), \partial B_i = \xi^{-1}(\{z_i\})$. The function $\xi$ has numerous names in the literature, including the *reaction co-ordinate*, and is typically defined such that higher values represent locations closer to $B$, as it is throughout this paper.

We take $(X(t) : t \geq 0)$ to be a continuous-time Markov process with almost surely continuous sample paths (the central example of diffusion processes will be the focus of our algorithmic constructions). We take $X(0)$ to be distributed according to an initial distribution $\lambda$, where the support of $\lambda$ is contained in $\xi^{-1}\left([z_A, z_B]\right)$ (so $X$ may begin on the boundary $\partial A$ of $A$; this choice is made for notational convenience, and the support of $\lambda$ may be taken instead to be all of $\mathbb{R}^d$ with only minor modifications). Except where necessary, dependence upon $\lambda$ will be suppressed from our notation.

Write $\tau_i$ for $\tau_{B_i}$, and define $\sigma_i = \tau_A \wedge \tau_i$ to be the first hitting time of $A \cup B_i$ for $X$. Note that $\tau_A, \tau_i, \sigma_i$ are equivalently the hitting times for the one-dimensional process $[\xi(X)]_s = \xi(X_s)$ of $\{z_i\}, \{z_A\}$ and $\{z_A, z_i\}$ respectively. In an algorithmic implementation, it is more convenient to use $\xi(X)$ to decide when a level has been crossed if $X$ has dimension greater than one. We remark that the process $\xi(X)$ is not in general a Markov process, so it is not possible to reduce all problems of this form to the univariate Markovian process setting.

The idea of MLS is to run a particle system in which each particle splits into several, say $R_i$, copies immediately upon first reaching a set boundary $\partial B_i$. Alternatively, it is terminated upon reaching $\partial A$ (See Fig. 1). The splitting is managed so that a healthy system of particles is available for estimating each $p_i$. Since $p = \prod_{i=1}^{m} p_i$, if we have an estimator $\hat{p}_i$ for each $p_i$, then $\prod_{i=1}^{m} \hat{p}_i$ is a natural estimator for $p$. The $\hat{p}_i$ are defined as follows: suppose that we begin with $N_0$ particles, of which $N_1$ reach level $B_1$ before $A$. Then we estimate $p_1$ with

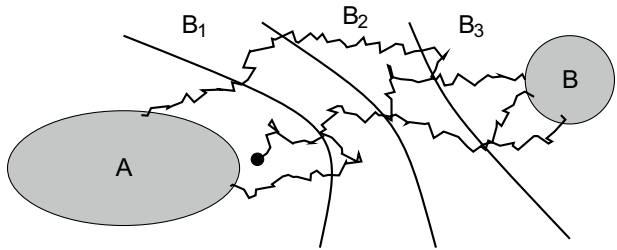

**Fig. 1** An illustration of multilevel splitting for a single particle system. The particle begins at the black node, and each branch splits into two i.i.d copies upon reaching a set $B_i$ for the first time. Branches which reach $A$ terminate there. Those which reach $B$ are used to form an estimate of the rare event probability $\mathbb{P}$ (process hits $B$ before A)

$$\hat{p}_1 = \frac{N_1}{N_0}.$$

Suppose that each of these $N_1$ surviving particles is split into a constant number $R_1$ of copies immediately upon reaching $B_1$, and that $N_2$ of the branched trajectories reach $B_2$ before $A$. Then we estimate $p_2$ with

$$\hat{p}_2 = \frac{N_2}{R_1 N_1}.$$

Continuing in this way, one obtains a sequence of estimators $\hat{p}_1, \ldots, \hat{p}_m$ of $p_1, \ldots, p_m$ which may be multiplied together to give the estimate

$$\hat{p} = \prod_{i=1}^{m} \hat{p}_i = \frac{N_m}{N_0 \prod_{i=1}^{m-1} R_i}$$

for $p$.

This estimator is unbiased; a proof may be found in [Section 3, Proposition 3.1] Amrein and Künsch (2011). We use a similar argument in Proposition 3.3 in Sect. 3.3 to establish the unbiasedness of the splitting-type algorithm which we develop there.

Provided that $z_i$ and $R_i$ are well-chosen, this estimate can be much more efficient than naïve Monte Carlo estimation. For instance, choosing $z_i$ such that $p_i = p^{1/m}$, and choosing $R_i$ with small variance such that $\mathbb{E}[R_i] = p_i^{-1}$, the relative variance of $\hat{p}$ is reduced to approximately $\mathcal{O}(p^{-1/m})$. See Glasserman et al. (1999); Lagnoux-Renaudie (2006), for more detail on parameter choice and asymptotic variance calculations.

It is convenient to work with the discrete-time pair-process $U_i = \big(\sigma_i, X_i(\sigma_i)\big)$, $i = 1, \ldots, m$, i.e. the values of X at its splitting times. Let $\mathcal{S}$ be the Borel sigma algebra associated with $\mathbb{R}_{\geq 0} \times \mathbb{R}^d$, and let $M_i : \big(\mathbb{R}_{\geq 0} \times \mathbb{R}^d\big) \times \mathcal{S} \to [0,1]$ denote the Markov kernels of this discrete-time process. Finally, we define potential functions $G_i : \mathbb{R}_{\geq 0} \times \mathbb{R}^d \to \{0,1\}$ on the state space of this pair process as indicators of the sets $B_i$ for the process $(U_i)$:

$$G_i(t,x) = \begin{cases} 1, & \text{if } \xi(x) \geq z_i, \\ 0, & \text{otherwise.} \end{cases}$$

A full description of multilevel splitting is given in Algorithm 1.

---

**Algorithm 1** Idealised Multilevel Splitting

Given $\lambda$ together with $G_i, M_i$ for $i = 1, \ldots, m$, an initial number of particles $N_0$, and splitting ratios $R_1, \ldots, R_{m-1}$:

1. For each $j = 1, \ldots, N_0$, draw independently $X_1^j(0) \sim \lambda$ and $U_1^j \sim M_1 \left( \big(0, X_1^j(0)\big), \cdot \right)$.

2. Let $S_1 = \{U_1^j : G_1(U_1^j) = 1\}$ be a list of the surviving paths, and set $N_1 = |S_1|$.

3. For $i = 2, \ldots, m$:

   (a) If $N_{i-1} = 0$, return $\hat{p} = 0$.

   (b) Else given $S_{i-1} = \{\tilde{U}_{i-1}^j\}_{j=1}^{N_{i-1}}$, for all $(j,k) \in \{(j',k') : 1 \leq j' \leq N_{i-1}, 1 \leq k' \leq R_{i-1}\}$ sample independently $U_i^{j,k} \sim M_i(\tilde{U}_{i-1}^j, \cdot)$.

   (c) Let $S_i = \{U_i^{j,k} : G_i(U_i^{j,k}) = 1\}$, and $N_i = |S_i|$.

4. Return

$$\hat{p} = \frac{N_m}{N_0 \prod_{i=1}^{m-1} R_i}.$$

---

A small modification to Algorithm 1 gives a variant commonly known as *fixed effort splitting*, which can be viewed as a Sequential Monte Carlo algorithm. This connection has been exploited previously by Cérou et al. (2006) (note that these algorithms are distinct from those which use SMC to approximate static rare events which depend upon the trajectory of a process only over a fixed time interval - see Cérou et al. (2012); Del Moral and Garnier (2005); Johansen et al. (2006)). In this variant, a number of particles $N$ to be maintained throughout is chosen in advance, and the splitting of each individual surviving particle is replaced with resampling from the set of surviving particles. This is useful in that the procedure does not require the specification of tuning parameters $R_1, \ldots, R_{m-1}$. A full description is given in Algorithm 2.

---

**Algorithm 2** Idealised MLS via SMC

Given $\lambda$ together with $G_i, M_i$ for $i = 1, \ldots, m$, and a fixed number of particles $N$:

1. For each $j = 1, \ldots, N$, draw independently $X_1^j(0) \sim \lambda$ and $U_1^j \sim M_1\left(\left(0, X_1^j(0)\right), \cdot\right)$.

2. Record $N_1 = \sum_{j=1}^N G_1(U_1^j)$.

3. For $i = 2, \ldots, m$:

   (a) If $N_{i-1} = 0$, return $\hat{p}^{\mathrm{SMC}} = 0$.

   (b) For $j = 1, \ldots, N$ sample independently $U_i^j \sim \frac{1}{N_{i-1}} \sum_{k=1}^N G_{i-1}(U_{i-1}^k) M_i(U_{i-1}^k, \cdot)$

   (c) Record $N_i = \sum_{j=1}^N G_i(U_i^j)$.

4. Return

$$\hat{p}^{\mathrm{SMC}} = \prod_{i=1}^m \left(\frac{N_i}{N}\right).$$

---

The particle system in Algorithm 2 has the familiar structure of an SMC sampler, with step 3b combining multinomial resampling with propagation via $M_i$. The numerical example in Sect. 4.2 is carried out using this simple scheme, but it should be noted that other resampling schemes from the SMC literature can also be used; see Gerber et al. (2019) for a detailed analysis of schemes which might be expected to reduce estimator variance without introducing any bias. The rare event probability of interest can be interpreted as the normalizing constant of an excursion-valued Feynman-Kac flow in the sense of [Section 12.2.6] Del Moral (2004). This flow has transition densities specified in terms of the underlying dynamics and stopping times, and zero-one-valued potential functions indicate whether crossing occurs into $B_i$ or $A$ at each level. Consequently, the SMC variant of MLS admits an interpretation as a mean field approximation of this flow and the estimator benefits from the usual theoretical analysis of these, see Del Moral (2004). This includes inheriting a strong law of large numbers, a central limit theorem and a proof of unbiasedness. This theory does not apply directly, however, to the estimator of Algorithm 1.

The unbiasedness proof of Amrein and Künsch (2011) also applies to Algorithm 2. An alternative but more general point of view from which this derives is the general theory of SMC estimators in the Feynman-Kac framework described above; see in particular [Theorem 7.4.2] Del Moral (2004).

**Remark 1** We have presented all of the algorithms in this paper from the perspective of estimating rare event probabilities. One may be interested in estimating other quantities, such as the law of the process conditional on its hitting $B$ before $A$. As with standard MLS, these may be estimated by a direct extension of the splitting algorithm.

## 2.2 Discretisation Error in MLS Algorithms

So far, we have assumed that we are able to simulate without approximation the pair $U_i = (\sigma_i, X(\sigma_i))$. Since these depend upon full sample paths of $X$, it is not apparent that this can be done except when $X$ has an exceptionally simple form, for example $X$ is a piece-wise deterministic process. In practice it is usual to resort to a discretisation scheme. For example, suppose that $X$ is described by $X(0) \sim \lambda$ and

$$dX(t) = \mu(X(t))dt + \sigma(X(t))dW(t) \tag{1}$$

for $t \in [0, T]$, where $W$ is $e$-dimensional Brownian motion for some $e \in \mathbb{N}$, and $\mu : \mathbb{R}^d \to \mathbb{R}^d, \sigma : \mathbb{R}^d \to \mathbb{R}^{d \times e}$ are sufficiently regular to guarantee the existence of a strong Itô solution (see for example, [Section 4.5] Kloeden and Platen (2013) for suitable conditions).

We might use an Euler-Maruyama scheme such as the following, defined on a chosen time-grid $t_j \in \mathscr{P}$ for a partition $\mathscr{P}$ of $[0, T]$:

$$\hat{X}(t_{j+1}) = \hat{X}(t_j) + \mu(\hat{X}(t_j))(t_{j+1} - t_j) + \sigma(\hat{X}(t_j))(W(t_{j+1}) - W(t_j)). \tag{2}$$

with $\hat{X}(0) \sim \lambda$. Such a scheme can then be used to implement an approximation of Algorithm 1 as follows: rather than drawing samples from $M_i$ in Steps 1 and 3b, one runs the discrete scheme until a crossing into $A$ or $B$ is observed at time $\hat{\sigma}_i = \min_j\{t_j : \hat{X}(t_j) \in A \cup B\}$, and approximates $U_i$ using $(\hat{\sigma}_i, \hat{X}(\hat{\sigma}_i))$.

This is the case for instance in Cérou and Guyader (2007), in which a formal algorithm is developed in a continuous-time setting but the numerical example is discretised "[finely] enough to avoid clipping the process, which could introduce a bias in the estimation". Lagnoux-Renaudie (2009) acknowledges explicitly the bias induced by discretisation in their application, and proposes a small modification to reduce, but not eliminate, it. Even Bréhier et al. (2016), which focuses on establishing the unbiasedness of a particular adaptive multilevel splitting framework in some generality ultimately invokes time-discretisation to apply the framework to continuous-time processes such as over-damped Langevin diffusions.

With such a numerical scheme, one is forced to assess the level-crossing problem according to the discrete sample paths of $\hat{X}$. But the law $\tilde{P}$ of $\hat{X}$ will not in general coincide with the true finite-dimensional marginal law $P$ of $(X_{t_1}, \dots, X_{t_k})$ induced by (1). And even if it were possible to get a finite-dimensional sample from $P$ restricted to the times of the partition, for example by exact simulation, this would give no information about the sample path over the open intervals $(t_j, t_{j+1})$, during which a crossing may (or may not) occur.

The quantities that are needed to carry out Algorithms 1 and 2 are $U_i^j = \left(\sigma_i^j, X_i^j\left(\sigma_i^j\right)\right)$ and $G_i(U_i^j)$. Using $\varepsilon$-strong simulation, we show that it is possible to sample exactly $G_i(U_i^j)$ without access to $U_i^j$ itself, using instead an $\varepsilon$-strong approximation to it. This is the focus of the next section; later, we show also that the modification to Algorithms 1 and 2 made necessary by using this approximation does not affect the unbiasedness of the resulting estimates.

## 2.3 $\varepsilon$-strong Simulation

Formally, following the definition given in Blanchet et al. (2017), an $\varepsilon$-strong algorithm is a joint construction of $X$ together with a family of processes $\tilde{X}^\epsilon$ indexed by $\epsilon > 0$ (defined on the same probability space) over an interval $[s, t]$ such that the following four properties hold:

1. Almost surely, $\sup_{r \in [s,t]} \|X(r) - \tilde{X}^\epsilon(r)\| \leq \varepsilon$ for an appropriate norm $\| \cdot \|$;
2. $\tilde{X}^\epsilon$ is piece-wise constant and left-continuous on $[s, t]$, taking only finitely many values and so can be fully stored on a computer;
3. $\tilde{X}^\epsilon$ can be simulated exactly. That is, to sample $\tilde{X}^\epsilon$ it is necessary to sample certain intermediate random variables, and this criterion requires that this can be done without approximations; and
4. Given a finite sequence of tolerances $\varepsilon_1 > \varepsilon_2 > \cdots > \varepsilon_m > 0$, for $1 \leq \ell_1 < \ell_2 \leq m$ it holds almost surely for all $r \in [s, t]$ that

$$\{x : \|\tilde{X}^{\epsilon_{\ell_2}}(r) - x\| \leq \varepsilon_{\ell_2}\} \subset \{x : \|\tilde{X}^{\epsilon_{\ell_1}}(r) - x\| \leq \varepsilon_{\ell_1}\},$$

and moreover it is possible to sample explicitly $\tilde{X}^{\epsilon_{\ell_2}}$ conditional on $\tilde{X}^{\epsilon_{\ell_1}}$.

An $\varepsilon$-strong algorithm produces a chain (in time) of finitely many $\| \cdot \|$-balls, each of which almost surely constrains the sample path of $X$ over the corresponding interval of time. Moreover, by applying Property 4 the radius of these balls can be iterative reduced, constraining $X$ progressively more tightly by employing a greater number of balls. An example (in two spatial dimensions) of how $\varepsilon$-strong sampling may be used is given in Fig. 2. It is often advantageous to apply condition 4 selectively in order to get tight constraints on $X$ at certain locations of interest, while allowing looser constraints elsewhere. For example, Fig. 2(c) shows the result of applying Property 4 to the first two $\varepsilon_1$-balls of the initial $\varepsilon_1$-sample in Fig. 2(b).

**Remark 2** The choice of a weak inequality in 1) differs slightly from the presentation in Blanchet et al. (2017). The reason is simply that our application requires calculating suprema and infima of the continuous function $\xi$ over regions $C(t) = \{x : \|\bar{X}(t) - x\| \leq \varepsilon\}$, and the weak inequality ensures that these extrema are attained in the regions $C(t)$.

**Remark 3** The insistence on condition 2) that $\tilde{X}^\epsilon$ be piece-wise constant is not strictly necessary since other processes which admit finite-dimensional representations could fill the same role. For example, continuous and piece-wise linear/polynomial $\tilde{X}^\epsilon$ are possible alternatives. However, we will assume throughout that $\tilde{X}^\epsilon$ both for convenience and for consistency with the $\varepsilon$-strong schemes which are employed in Sect. 4.

**Remark 4** In existing multi-dimensional $\varepsilon$-strong algorithms, the norm in condition 1) is always the supremum norm, so the representation in Fig. 2 corresponding to the Euclidean norm should be taken to be schematic.

With Property 4 in mind, let $(\varepsilon_\ell)_{\ell=1}^\infty$ be a decreasing sequence of tolerances converging to 0. Write $\tilde{X}^\ell[s : t]$ for the $\varepsilon$-strong path $\tilde{X}^{\epsilon_\ell}[s : t]$. Let $(t_1, \ldots, t_K)$ be the jump-times of this path, and let $t_0 = s$, $t_{K+1} = t$. It is useful to define the associated discrete-time process $\left(\tilde{X}_k^\ell\right)_{k=0}^K$ where $\tilde{X}_k^\ell = \tilde{X}^\ell(t_k)$. We define also an augmented process called the *skeleton* of $\tilde{X}^\ell$ to be the discrete-time process $\left(\tilde{Z}_k^\ell\right)_{k=0}^K$ such that

$$\tilde{Z}_k^\ell = \left(t_k, t_{k+1}, \tilde{X}_k^\ell, \ell\right).$$

(This is, in a way, rather a backwards definition since an $\varepsilon$-strong path itself is typically constructed from its skeleton.) Given the skeleton $Z$ as defined above, the (almost) unique $\varepsilon$-strong path associated with it is defined by

$$\tilde{X}^\ell(u) = \sum_{k=0}^{K} \mathbb{I}_{[t_k, t_{k+1})}(u)\tilde{X}_k^\ell.$$

for $u \in [s, t)$, and we may take $\tilde{X}^\ell(t) = \tilde{X}_K^\ell$. The skeleton is somehow a more computationally-motivated object than its associated path, and we will refer primarily to the paths themselves outside of our algorithmic pseudo-code. It is useful to define also $C_k = \{x\ :\ \|x - \tilde{X}_k^\ell\| \leq \varepsilon_\ell\}$, the constraining region for $X$ over $[t_k, t_{k+1}]$.

Say two skeletons $\left(\tilde{Z}_k^{i,1}\right)_{k=0}^{K}$, $\left(\tilde{Z}_m^{j,2}\right)_{m=0}^{L}$ defined on $[r, s]$ and $[s, t]$ respectively are *compatible* if $C_K^1 \cap C_0^2 \neq \emptyset$. For two compatible skeletons, we define their concatenation $\tilde{Z}^3 = \tilde{Z}^{i,1} \oplus \tilde{Z}^{j,2}$ to be the process $\left(\tilde{Z}_n^3\right)_{n=1}^{K+L+1}$ with

$$\tilde{Z}_n^3 = \begin{cases} \tilde{Z}_n^{i,1}, & 0 \leq n \leq K, \\ \tilde{Z}_{n-K-1}^{j,2}, & K+1 \leq n \leq L+K+1. \end{cases}$$

Analogously, for two compatible $\varepsilon$-strong paths $\tilde{X}^{i,1}, \tilde{X}^{j,2}$ defined on $[r, s]$, $[s, t]$ respectively, and with skeletons $\left(\tilde{Z}_k^{i,1}\right)_{k=0}^{K}$, $\left(\tilde{Z}_m^{j,2}\right)_{m=0}^{L}$, we define their concatenation as follows. Take the skeleton $\tilde{Z}^3 = \tilde{Z}^{i,1} \oplus \tilde{Z}^{j,2}$, and writing $\tilde{Z}_n^3 = \left(t_n, t_{n+1}, \tilde{X}_n^{\ell(n)}, \ell(n)\right)$ set

$$\left(\tilde{X}^{i,1} \oplus \tilde{X}^{j,2}\right)(u) = \sum_{n=0}^{K+L+1} \mathbb{I}_{[t_n, t_{n+1})}(u)\tilde{X}_n^{\ell(n)}$$

for $u \in [s, t)$, $\left(\tilde{X}^{i,1} \oplus \tilde{X}^{j,2}\right)(t) = \tilde{X}_{L+K+1}^{\ell(L+K+1)}$.

Two paths which are themselves concatenations of $\varepsilon$-strong paths may be concatenated analogously. We will exploit the fact that the binary concatenation operation is associative to allow us to write concatenations of more than two processes without ambiguity. We also find it convenient to adopt the convention that for $k_1 \geq k_2$, the sub-skeleton $\left(\tilde{Z}_k\right)_{k=k_1}^{k_2}$ acts as the identity for this binary operation, so that for any skeleton $\tilde{Y}$, $(Z_k)_{k=k_1}^{k_2} \oplus \tilde{Y} = \tilde{Y} \oplus (Z_k)_{k=k_1}^{k_2} = \tilde{Y}$.

**Remark 5** Typically $\varepsilon$-strong algorithms require more information about the process than we have made explicit in our definition of a skeleton, for example those of Pollock et al. (2016). For ease of exposition, we have suppressed this since we do not need to refer to it for the development of the algorithms in this paper, but it should be understood that our skeletons contain any extra information required for the stated $\varepsilon$-strong conditions to hold.

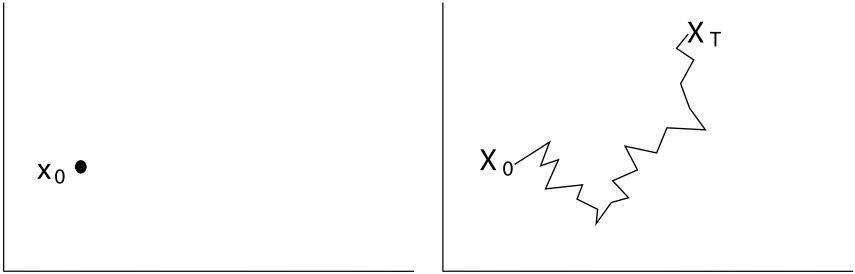

**(a)** Initial value $X_0 = x_0$, and target path $X$

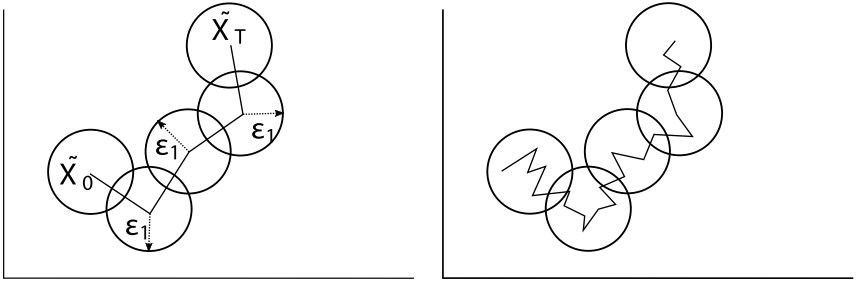

**(b)** $\varepsilon_1$-strong sample and constraining regions

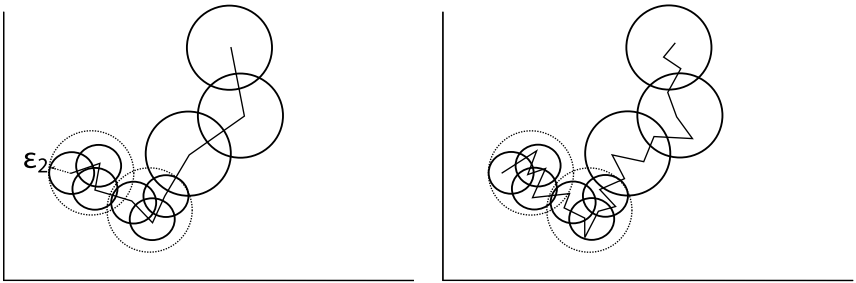

**(c)** Partial $\varepsilon_2$-strong sample

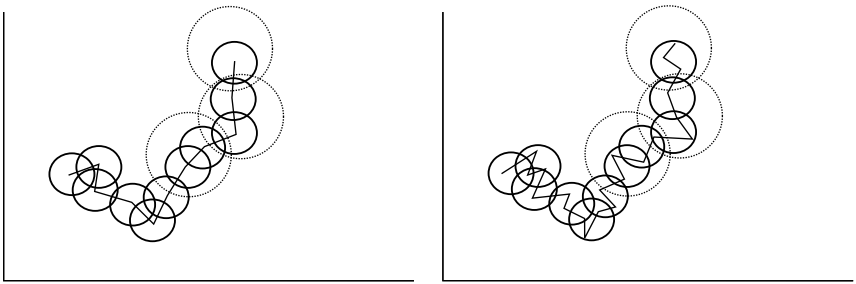

**(d)** Full $\varepsilon_2$-strong sample

## 3 Exact Simulation of Rare Events

In multilevel splitting, one tracks the progress of $X$ towards $A$ and $B$ using the reaction co-ordinate $\xi$, declaring a crossing at level $i$ when the process $\xi(X)$ reaches either $z_A$ or $z_i$. In this section we describe methods for sampling such *barrier crossing* events exactly, for Markov processes $X$ with almost surely continuous sample paths for which an $\varepsilon$-strong method for sampling $X$ exists. Combining these with a slight modification of Algorithms 1 & 2 provides a method of obtaining unbiased estimates of rare event probabilities. Much of what follows is geometrically intuitive, though notationally cumbersome, and Figs. 2 and 3 are intended to illustrate the intuition which motivates the accompanying specifications.

Throughout this section we take $X$ to be a diffusion as described in (1), together with the conditions assumed there. For simplicity, we assume further that $X$ has volatility bounded away from 0, which ensures that $X$ crosses any given boundary with positive probability over any time interval. We assume also that an $\varepsilon$-strong algorithm as described in Sect. 2.3 has been chosen and is used to carry out the sampling in Algorithms 3, 4 and 6.

### 3.1 Crossing a Single Barrier

We begin with the simpler problem of sampling exactly an indicator random variable for the event that $X$ crosses into a set $D = \xi^{-1}([z_D, \infty))$ when started from its complement $D^c = \xi^{-1}((-\infty, z_D))$, over the fixed time interval $[0, t]$. To this end, we suppose that $X(0) \sim \lambda$ where the support of $\lambda$ is contained in $D^c$. Assume that X is sufficiently regular that, almost surely, a path $X[0 : t]$ which crosses into $D$ attains a maximum distance $d^{\max}(X, D^c) > 0$ from $D^c$, and conversely, a path $X[0 : t]$ which does not cross into $D$ has (almost surely) minimum distance $d^{\min}(X, D) > 0$ from $D$. Consider the $\varepsilon$-strong path $\tilde{X}^{\varepsilon_1}(0 : t)$ for a tolerance $\varepsilon_1$, and let $0 = t_0 < \cdots < t_{K+1} = t$ be its jump-times. For $k = 0, \ldots, K$, inside each time interval $[t_k, t_{k+1}]$ the ball $C_k = \{x : \|x - \tilde{X}^{\varepsilon_1}(t_k)\| \leq \varepsilon_1\}$ almost surely constrains the path of $X$ associated with $\tilde{X}^{\varepsilon_1}$. So if $\varepsilon_1 < \max(d^{\max}(X, D^c), d^{\min}(X, D))$, then either i) $C_k \subset D$ for some $k$ (if $X$ does make a crossing), or ii) $C_k \subset D^c$ for all $k$ (if $X$ does not make a crossing). By checking each $C_k$ in turn, we can determine which of these conditions holds, and thereby construct the desired indicator random variable.

Of course, it is not possible to choose a suitable $\varepsilon_1$ in advance, since the underlying path $X$ and its minimum and maximum distances from $D$ and $D^c$ are not known. Instead, we can specify a sequence $(\varepsilon_\ell)_{\ell=1}^{\infty}$ of tolerances with $\varepsilon_\ell \to 0$. If $X^\ell$ turns out to be insufficient to determine the crossing, we can apply Property 4 of Sect. 2.3 to sample $X^{\ell+1}$ conditional on $X^\ell$ as necessary until a sufficiently small tolerance is found.

It is very wasteful, however, to construct a finitely-representable path $\tilde{X}[0 : t]$ which is very close to $X[0, t]$ on the whole interval $[0 : t]$. It is likely that even when $X$ crosses into $D$, much of the time $X$ is not near the boundary $\partial D$, and we need only approximate $X$ closely where it is near $\partial D$. For this reason, as suggested in Sect. 2.3, it useful to work instead with paths of mixed tolerance

$$\tilde{X} = \bigoplus_{j=1}^{J} \tilde{X}^{\ell(j)}[s_{j-1}, s_j],$$

where $(0 = s_0, s_1, \ldots, s_J)$ is a partition of $[0, t]$ with $s_J = t$, and $\left(\varepsilon_{\ell(j)}\right)_{j=1}^{J}$ is a selection from $(\varepsilon_\ell)_{\ell=1}^{\infty}$. Such a path is the result, for example, of applying Property 4 with $\varepsilon_2 < \varepsilon_1$ to a constant segment $\tilde{X}^1(t_{k-1}, t_k)$ of $\tilde{X}^1$, and the result in this case would be

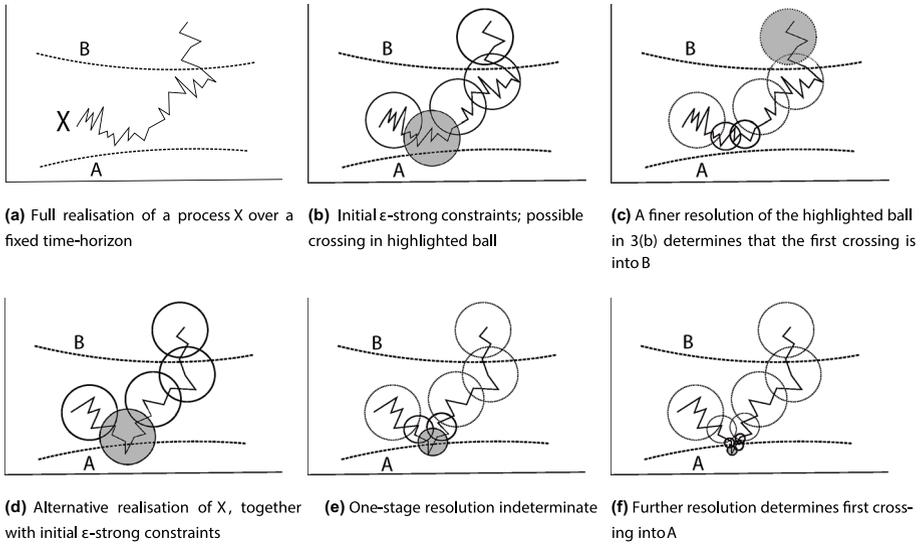

**(a)** Full realisation of a process X over a fixed time-horizon

**(b)** Initial ε-strong constraints; possible crossing in highlighted ball

**(c)** A finer resolution of the highlighted ball in 3(b) determines that the first crossing is into B

**(d)** Alternative realisation of X, together with initial ε-strong constraints

**(e)** One-stage resolution indeterminate

**(f)** Further resolution determines first crossing into A

**Fig. 3** Two illustrations of Algorithm 4. The first row shows: (**a**) a realisation of $X$ over a finite time horizon, (**b**) an initial $\varepsilon$-strong simulation and (**c**) a refinement which is sufficient to show the process crossing into $B$. The second row shows (**d**) an alternative sample path consistent with the same initial $\varepsilon$-strong simulation, (**e**) an inconclusive refinement and (**f**) a further refinement sufficient to conclude that the process has crossed into $A$

$$\tilde{X} = \tilde{X}^1[0 : t_{k-1}] \oplus \tilde{X}^2[t_{k-1} : t_k] \oplus \tilde{X}^1[t_k : t].$$

For later convenience, our formalisation in Algorithm 3 of the algorithm under description takes an $\tilde{X}$ of this kind, or rather the skeleton of such a path, as input.

The three possible relationships between $C_k$ and $D, D^c$ can be described in terms of $\xi$, the reaction co-ordinate: $\sup_{x \in C_k} \xi(x) < z_D$ is equivalent to $X[t_k, t_{k+1}] \subset D^c$, and similarly $\inf_{x \in C_k} \xi(x) \geq z_D$ is equivalent to $X[t_k, t_{k+1}] \subset D$. The third possibility, that

$$\inf_{x \in C_k} \xi(x) < z_D, \ \sup_{x \in C_k} \xi(x) \geq z_D \,,$$

gives no definite information about the location of $X[t_k, t_{k+1}]$ with respect to $D$. It is consistent with $X[t_k, t_{k+1}]$ falling entirely in $D$, entirely in $D^c$, or partially in both. We categorise the behaviour of the process in this time interval by defining:

$$n_k := \begin{cases} -1, & \text{if } \sup_{x \in C_k} \xi(x) < z_D \\ 0, & \text{if } \inf_{x \in C_k} \xi(x) < z_D, \ \sup_{x \in C_k} \xi(x) \geq z_D \\ +1, & \text{if } \inf_{x \in C_k} \xi(x) \geq z_D \end{cases} \tag{3}$$

Algorithm 3 samples exactly an indicator for the event that $X$ crosses into $D$.

---

**Algorithm 3** Single barrier crossing

function($(\tilde{Z}_k)_{k=0}^K, D$):

1. Calculate the sequence $(n_k)_{k=0}^K$.

2. If $n_k = -1$ for all $k = 1, \ldots, K$, return $(0, \tilde{Z})$ to indicate no crossing into $D$.

3. If $n_k = +1$ for some $k$, return $(+1, \tilde{Z})$ to indicate a crossing into $D$.

4. Else:

   (a) Set $j = \min\{k \in \{0, \ldots, K\} : n_k = 0\}$, and consider $\mathcal{Z}^{\ell(j)} := \tilde{Z}_j^{(j)} = (t_j, t_{j+1}, \tilde{X}_j, \ell(j))$. Use the refining Property 4 of Section 2.3 to sample $\mathcal{Z}^{\ell(j)+1}$ conditional on $\mathcal{Z}^{\ell(j)}$.

   (b) Update

   $$\tilde{Z} \leftarrow (\tilde{Z}_k)_{k=0}^{j-1} \oplus \mathcal{Z}^{\ell(j)+1} \oplus (\tilde{Z}_m)_{m=j+1}^K,$$

   and update $K \leftarrow (\#\text{jump-times of } \tilde{Z}) + 1$. Return to Step 1.

---

**Remark 6** It may be noted that in Step 4a) of Algorithm 3, it is not strictly necessary to choose $k$ minimal. There may be computational advantages to using a different system, such as attempting choose an $k$ for which $C_k \cap A$ is large (indicating a high probability of crossing). This can be computationally preferable, at the expense of providing less information about $\tau_A$ (see Sect. 3.3).

The assumption that $\sup_{x \in C_k} \xi(x)$ and $\inf_{x \in C_k} \xi(x)$ can be calculated is rather strong, but holds for many realistic scenarios. For example, supposing $X$ takes values in $\mathbb{R}^d$ and, taking the norm, $\|x\| = \max_{i \in 1, \ldots, d} |x_i|$, these quantities can be calculated if $\xi$ is monotonic in each argument. As a specific example, in Sect. 4, we take $d = 2$, $\xi(x, y) = \min(x, y)$. Another example of a tractable reaction coordinate, which illustrates that monotonicity is not necessary, is $\xi(x, y) = |x - y|$.

## 3.2 Crossing a Two-sided Barrier

We consider next a two-sided barrier problem, with regions $A = \xi^{-1}((-\infty, z_A])$, $B = \xi^{-1}([z_B, \infty))$, with $X(0) \sim \lambda$ such that $z_A \leq \xi(X(0)) < z_B$, and the problem of sampling an indicator random variable for the event that $X$ crosses into $B$ before $A$. Here we work over over a random interval $[0, \sigma]$ where $\sigma$ is the hitting time for $A \cup B$ of $X$, rather than over a fixed interval as in the previous section.

In this case, we can declare a level crossing into $A$ (for example) at the first $k$ for which

$$\sup_{x \in C_k} \xi(x) < z_A \text{ and } \max_{j < k} \sup_{x \in C_k} \xi(x) < z_B,$$

if such an $k$ exists.

Informally, we can declare the crossing into $A$ when i) some $\varepsilon$-ball lies entirely in set $A$, which guarantees that $X$ has reached $A$; and ii) no preceding $\varepsilon$-ball intersects set $B$, which guarantees that $X$ has not reached $B$. (The conditions for a crossing into $B$ are analogous).

**Table 1** Enumeration of possible $C_k$-locations with respect to a two-sided barrier

| $n_k$ | Condition | Meaning |
|---|---|---|
| −2 | $\sup_{x \in C_k} \xi(x) \leq z_A$ | $X$ remains within $A$ on $[t_k, t_{k+1}]$ |
| −1 | $\inf_{x \in C_k} \xi(x) \leq z_A$, $\sup_{x \in C_k} \xi(x) > z_A$ | $X$ may enter $A$ in $[t_k, t_{k+1}]$ |
| 0 | $\inf_{x \in C_k} \xi(x) > z_A$, $\sup_{x \in C_k} \xi(x) < z_B$ | $X$ does not enter $A$ or $B$ in $[t_k, t_{k+1}]$ |
| 1 | $\sup_{x \in C_k} \xi(x) \geq z_B$, $\inf_{x \in C_k} \xi(x) < z_B$ | $X$ may enter $B$ in on $[t_k, t_{k+1}]$ |
| 2 | $\inf_{x \in C_k} \xi(x) \geq z_B$ | $X$ remains within $B$ on $[t_k, t_{k+1}]$ |

If there is no such $k$, it is necessary to carry out further simulations using Property 4 of Sect. 2.3.

As in the previous section, we associate a number $n_k \in \{-2, -1, 0, 1, 2\}$ with each ball $C_k$, according to the categorisation in Table 1. In order to simplify the categorisation and presentation of the algorithm, we make the assumption that our initial tolerance $\varepsilon_1$ is sufficiently small that i) $C_k$ intersects at most one of $A$, $B$; and that ii) if $C_k \cap A \neq \emptyset$, then $C_{k+1} \cap B = \emptyset$ (likewise with $A$, $B$ interchanged); this can be assured by refining the initial tolerance until it is satisfied. With this assumption, the categorisation in Table 1 is complete, and the sequence $(n_k)$ satisfies $n_{k+1} \in \{n_k - 1, n_k, n_{k+1}\}$.

Suppose we have calculated the sequence $(n_k)_{k=0}^K$ associated with the path $\tilde{X}[s:t]$. In order to determine which of $A$ and $B$ has been crossed first, it is necessary to consider segments of $\tilde{X}[s:t]$ in which a crossing of one or the other barrier may have occurred, but a crossing of both barriers *cannot* have occurred. By checking each of these segments in turn, the decision can be made. In terms of the sequence $(n_k)$, these segments are constructed as follows. We write $J$ for the number of segments, where the definition of $J$ is contained in the construction. We define recursively the sequence of indices which mark the beginning of a new segment in which a crossing may occur, as the sequence $(\kappa(j))_{j=0}^J$. Let $\kappa(0) = 0$, and while $\kappa(j-1) < K + 1$, set

$$\kappa(j) = \begin{cases} \min\{k > \kappa(j-1) : n_k = 0\} \wedge (K+1) \text{ if } n_{\kappa(j-1)} \neq 0 \\ \min\{k > \kappa(j-1) : n_k \neq 0\} \wedge (K+1) \text{ if } n_{\kappa(j-1)} = 0. \end{cases}$$

Each element of this sequence is taken to denote the beginning of a block $\mathcal{B}_j \subset (n_k)$ of consecutive elements, so $\mathcal{B}_j = \{n_k : \kappa(j-1) \leq k < \kappa(j) - 1\}$. By construction, each $\mathcal{B}_j$ consists of a string of elements of exactly one of the sets $\{-2, -1\}, \{0\}, \{1, 2\}$. Each block therefore corresponds to a segment of $\tilde{X}[s:t]$ in which $X$ crosses into at most one of $A$ and $B$. For example, in the case that $(n_k) = (1, 1, 0, 0, -1, -2, -1)$, $J = 3$ and the blocks are $\mathcal{B}_1 = (1, 1), \mathcal{B}_2 = (0, 0), \mathcal{B}_3 = (-1, -2, -1)$, which in this case correspond to "possible crossing into $B$", "no crossing" and "definite crossing into $A$", respectively.

The two-sided barrier crossing procedure is given in Algorithm 4, in which the output is an indicator random variable for the event that $X$ hits $B$ before $A$. An illustration is given in Fig. 3.

---

**Algorithm 4** Two-sided barrier crossing

function$((\tilde{Z}_k)_{k=0}^K, A, B)$, for $\tilde{Z}$ a skeleton over the interval $[s,t]$:

1. Initialise $\tilde{Z}^{\text{full}}$ as an empty skeleton.

2. Calculate the sequence $(n_k)_{k=0}^K$ associated with $\tilde{Z}$.

3. Divide $(n_k)_{k=0}^K$ into blocks $\mathscr{B}_1, \ldots, \mathscr{B}_J$ using the sequence $(\kappa(j))_{j=0}^J$ as described above.

4. For $j = 1, \ldots, J$:

   (a) If $(-2) \in \mathscr{B}_j$, set $D = (-1)$ to indicate a crossing into $A$, and skip to 6.

   (b) If $(+2) \in \mathscr{B}_j$, set $D = (+1)$ to indicate a crossing into $B$, and skip to 6.

   (c) If $(-1) \in B_j$, sample $(I, \tilde{Z}') \leftarrow$ Algorithm.3 $((\tilde{Z}_k)_{k=\kappa(j-1)}^{\kappa(j)-1}, A)$ to decide if the first crossing is into $A$, and set

   $$\tilde{Z} \leftarrow (\tilde{Z}_k)_{k=0}^{\kappa(j-1)-1} \oplus \tilde{Z}' \oplus (\tilde{Z}_k)_{k=\kappa(j)}^K,$$

   and $K \leftarrow$ (#jump-times of $\tilde{Z}$) + 1. If $I = (+1)$, set $D = (-1)$ and skip to 6.

   (d) If $(+1) \in B_j$, sample $(I, \tilde{Z}') \leftarrow$ Algorithm.3 $((\tilde{Z}_k)_{k=\kappa(j-1)}^{\kappa(j)-1}, B)$ to decide if the first crossing is into $B$, and update $\tilde{Z}, K$ as in c). If $I = (+1)$, set $D = (+1)$ and skip to 6.

5. Here, we know that no crossing is made in the interval $[s,t]$ spanned by $\tilde{Z}$, so it is necessary to sample a continuation of the path. Let $\tilde{Z}^{\text{full}} \leftarrow \tilde{Z}^{\text{full}} \oplus \tilde{Z}$. Writing $(t_K, t, x, \varepsilon) = \tilde{Z}_K$, update $(s,t) \leftarrow (t, 2t - s)$, and sample a new $\tilde{X}^{\varepsilon_1}(s:t]$. Record its skeleton $(\tilde{Z}_k)_{k=1}^K$, and return to 2.

---

6. Set $\tilde{Z}^{\text{full}} \leftarrow \tilde{Z}^{\text{full}} \oplus \tilde{Z}$, and return $(D, \tilde{Z}^{\text{full}})$.

---

**Remark 7 (Implementation)** In general, it may be computationally inefficient to use a sufficiently small initial $\varepsilon$ for the whole sample path of $X$, for example when $X$ crosses a barrier at a very early time. We note that there are many variations of Algorithm 4 which will also sample the outcome correctly and may avoid doing so for computational efficiency. Our choice has been made for clarity of exposition.

**Remark 8 (Discontinuous processes and jump-diffusions)** We have assumed throughout for convenience that the sample paths with which we deal are almost surely continuous. Relaxing this requirement is straightforward but slightly complicates the implementation. Given an $\varepsilon$-strong algorithm for a jump-diffusion or similar piece-wise-continuous process, an appropriate alteration to the rule for beginning a new block will give an equally correct algorithm.

## 3.3 Exact MLS

Finally, we turn to an exact implementation of multilevel splitting. The main point of difference with Algorithm 1 is that since an $\varepsilon$-strong sample $\tilde{X}[s,t]$ merely constrains the corresponding path $X[s,t]$, there is no easy way to determine the hitting location and time of any given barrier. In particular, we will not have access to the hitting times $\sigma_i$, $\tau_i$ nor the hitting locations $X(\tau_i)$ defined in Sect. 2.1.

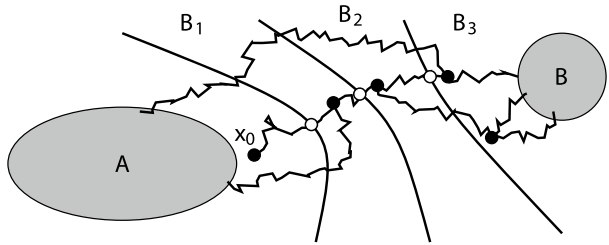

**Fig. 4** An illustration of Idealised Splitting with Couplings for a single particle system. The particle begins at the node labelled $x_0$. Level crossings are indicated by empty nodes, whereas splittings occur at the filled nodes. Between any empty node and the following filled node, the particle trajectories are coupled identically. Compare with Fig. 1

Suppose we use Algorithm 4 to sample an indicator random variable for the event that $X$ hits $B_1$ before $A$, for instance, with initial simulation interval $[0, T]$. Suppose that a positive result is returned over the interval $[0, cT]$, for some random $c \in \mathbb{N}$ corresponding to the number of passes through Algorithm 4. We must then choose when and where to split this path of $\tilde{X}$. In this section, we show that if the splitting is carried out at time $cT$, this does not affect the unbiasedness of the MLS estimate.

In general, write $\tilde{\sigma}_i$ for a random time which serves as an upper bound on the first hitting time of $A \cup B_i$ for $\tilde{X}$, which is defined as:

$$\tilde{\sigma}_i = T \cdot \min\{m \in \mathbb{N} : mT \geq \sigma_i\},$$

i.e. the time to which $\tilde{X}$ is sampled in Algorithm 4 (so $\tilde{\sigma}_i$ is a multiple of $T$) in order to establish that a crossing into $A$ or $B_i$ has occurred. Similarly, let $\tilde{\tau}_i$ be the corresponding upper bound on the first hitting time of $B_i$. To understand the exact MLS we describe later in this section, it is helpful to have in mind an idealised splitting scheme slightly different from MLS as presented in Sect. 2.1, which we call *idealised splitting with coupling*.

As in idealised MLS, we assume that it is possible to sample complete continuous paths of $X$ up to a given stopping time. But rather than split these paths into *independent* copies at times $\tau_i$, the split paths are coupled in the following way: from time $\tau_i$ until time $\tilde{\tau}_i$, the "split" paths are set to be identically equal, and after this time they evolve conditionally independently given $X_{\tilde{\tau}_i}$.

For $i = 1, \ldots, m$, let $\tilde{M}_i$ denote the transition kernels for the discrete time quadruple process $V_i = (\sigma_i, \tilde{\sigma}_i, X(\sigma_i), \tilde{X}(\tilde{\sigma}_i))$. Define also $\tilde{G}_i(V_i) = \mathbb{1}_{B_i}(X(\sigma_i))$. The details are given in Algorithm 5. Call the estimator for $p$ resulting from this algorithm $\tilde{p}$. (See Fig. 4)

Of course, it is not possible to implement Algorithm 5 as written, since we cannot simulate full paths of $X$, nor make splits at times $\tau_i$. But the construction of MLS with couplings means that an algorithm which splits paths at the tractable time $\tilde{\tau}_i$ instead (and allows them to propagate independently from that point) produces identical estimators $\tilde{q}_i$ for $p_i$.

It is possible that a particle crosses several barriers $z_i, z_{i+1}, \ldots$ before time $\tilde{\tau}_i$. In this case, the particle splits as normal at each barrier, each new copy remaining identically coupled, until the splitting time $\tilde{\tau}_i$ after which they proceed independently.

The following proposition establishes that this framework gives rise to unbiased estimates. In order to simplify the analysis, we make a further assumption about the $\varepsilon$-strong method being used: that over any interval $[s, t]$, we have $\left(\tilde{X}(t) \mid \tilde{X}(s) = x_s\right)$ is equal in distribution to $\left(X(t) \mid X(s) = x_s\right)$. In other words, we assume that the end-points of $\varepsilon$-strong samples are from the true distribution of $X$. This holds in the schemes developed in Beskos et al. (2012) and Pollock et al. (2016), for example. And given the close relation of $\varepsilon$-strong simulation with exact simulation, it may often be the case that a $\varepsilon$-strong algorithm of this form can be constructed (see for instance Blanchet et al. (2017); Blanchet and Zhang (2020)).

Given the MLS with couplings scheme, it is now possible to establish this result with only minor modification of the existing arguments of Amrein and Künsch (2011); this result paves the way for the methodology which we introduce and in principle allows for unbiased estimation of rare event probabilities in continuous time whenever $\varepsilon$-strong simulation of the process of interest is possible.

---

**Algorithm 5** Idealised Splitting with Coupling

Given $\lambda$ together with $\tilde{G}_i, \tilde{M}_i$ for $i = 1, \ldots, m$, an initial number of particles $N_0$, and splitting ratios $R_1, \ldots, R_{m-1}$:

1. For $j = 1, \ldots, N_0$:

   (a) Draw $X_1^j(0) \sim \lambda$, $V_1^j \sim \tilde{M}_1\left(\left(0, 0, X_1^j(0), X_1^j(0)\right), \cdot\right)$.

2. Let $S_1 = \{V_1^j : \tilde{G}_1(V_1^j) = 1\}$ be a list of the the surviving paths, and $N_1 = |S_1|$.

3. For $i = 2, \ldots, m$:

   (a) If $N_{i-1} = 0$, return $\tilde{p} = 0$.

   (b) Otherwise given $S_{i-1} = \{V_{i-1}^j\}_{j=1}^{N_{i-1}}$, for each $(j,k) \in \{(j',k') : 1 \leq j' \leq N_{i-1}, 1 \leq k' \leq R_{i-1}\}$:

      i. Sample $V_i^{(j,k)} \sim \tilde{M}_i(V_{i-1}^j, \cdot)$

   (c) Let $S_i = \{V_i^{j,k} : \tilde{G}_i(V_i^{j,k}) = 1\}$, and set $N_i = |S_i|$.

4. Estimate

$$\tilde{p} = \frac{N_m}{N_0 \prod_{i=1}^{m-1} R_i}.$$

---

**Proposition 1** *$\tilde{p}$ is an unbiased estimator for $p$: $\mathbb{E}[\tilde{p}] = p$.*

**Proof** We consider the discrete-time Markov process $(V_i)_{i=0}^m$. Let $\tilde{M}_i$ denote the transition kernels of this process at time $i$. Let also $\tilde{M}_{i:j} := \tilde{M}_j \circ \tilde{M}_{j-1} \circ \ldots \circ \tilde{M}_{i+1}$ denote the elements of the associated two-parameter dynamic semigroup, describing the evolution of the process from time $i$ to $j$. Note that for every $i$, $\Delta = \mathbb{R}^2 \times A \times \mathbb{R}^d$ is an absorbing state for $\tilde{M}_i$ since if $X(\sigma_i) \in A$, $\sigma_j = \sigma_i$ for all $j > i$.

Define $\mathcal{G}_0$ to be the sigma algebra generated by $\{V_0^j : j = 1, \ldots, N_0\}$, $\mathcal{G}_k = \mathcal{G}_{k-1} \vee \Sigma(V_k^j : 1 \leq j \leq N_k)$ for $k = 2, \ldots, m$, so that $(G_k)_{k=1}^m$ is the natural filtration of this process..

We observe the following recursion for the number of particles successfully reaching $B_i$ given $\mathcal{G}_{i-1}$, which is immediate from the definition of $M_i$: that for any function $h_k : \mathbb{R}^2 \times (\mathbb{R}^d)^2 \to \mathbb{R}$ which is equal to 0 on $\Delta$:

$$\mathbb{E}\left[\sum_{j=1}^{R_{i-1}N_{i-1}} h_i\left(V_i^j\right) \middle| \mathcal{G}_{i-1}\right] = R_{i-1} \sum_{k : \tilde{G}_{i-1}(V_{i-1}^k)=1} \int h_i(u) \tilde{M}_i\left(V_{i-1}^k, du\right). \tag{4}$$

(taking $R_0 = 1$) which follows since each $V_i^j$ which is a descendent of any particular $V_{i-1}^k$ has the same marginal law.

Define estimators $\tilde{p}_i$ for $p_i$ in Algorithm 5, namely

$$\tilde{p}_i = \begin{cases} \frac{N_1}{N_0}, & i = 1, \text{ and} \\[2ex] \frac{N_i}{R_{i-1}N_{i-1}}, & 2 \le i \le m. \end{cases}$$

We show that for $1 \le k$,

$$\mathbb{E}\left[\prod_{i=k}^m \tilde{p}_i \Big| \mathcal{G}_{k-1}\right] = \frac{1}{N_{k-1}} \sum_{j\,:\,\tilde{G}_{k-1}\left(V_{k-1}^j\right)=1} \left(1 - \tilde{M}_{(k-1):m}\left(V_{k-1}^j, \Delta\right)\right),$$

by backwards induction on $k$, starting with the case $k = m$. Note that for the case $k = 1$, each term on the sum on the RHS is the probability that a particle with a given starting value successfully reaches $B$ before $A$. The result is then obtained upon taking a further expectation over the starting value.

The case $k = m$:

$$\tilde{p}_m = \frac{N_m}{R_{m-1}N_{m-1}}$$

$$= \frac{1}{R_{m-1}N_{m-1}} \sum_{j=1}^{R_{m-1}N_{m-1}} \tilde{G}_m\left(V_m^j\right)$$

and the result then follows from taking the conditional expectation, and applying (4) with $h_j\left(\sigma, \tilde{\sigma}, X(\sigma), \tilde{X}(\tilde{\sigma})\right) = \mathbb{I}_{B_j}(X(\sigma))$.

Now supposing the result holds for $k+1$, we show that it holds also for $k$. We have the following chain of equalities (with the convention that $G_0\left(U_0^j\right) = 1$):

$$\mathbb{E}\left[\prod_{i=k}^m \tilde{p}_i \Big| \mathcal{G}_{k-1}\right] = \mathbb{E}\left[\tilde{p}_k \mathbb{E}\left[\prod_{i=k+1}^m \tilde{p}_i \Big| \mathcal{G}_k\right] \Big| \mathcal{G}_{k-1}\right] \tag{5}$$

$$= \mathbb{E}\left[\frac{N_k}{R_{k-1}N_{k-1}} \cdot \frac{1}{N_k} \sum_{j\,:\,\tilde{G}_k\left(V_k^j\right)=1} \left(1 - \tilde{M}_{k:m}\left(V_k^j, \Delta\right)\right) \Big| \mathcal{G}_{k-1}\right] \tag{6}$$

$$= \frac{1}{N_{k-1}} \sum_{j\,:\,\tilde{G}_{k-1}\left(V_{k-1}^j\right)=1} \left(\int \left(1 - \tilde{M}_{k:m}(u, \Delta)\right)\tilde{M}_k\left(V_{k-1}^j, du\right)\right) \tag{7}$$

$$= \frac{1}{N_{k-1}} \sum_{j\,:\,\tilde{G}_{k-1}\left(V_{k-1}^j\right)=1} \left(1 - \int \tilde{M}_{k:m}(u, \Delta)\tilde{M}_k\left(V_{k-1}^j, du\right)\right) \tag{8}$$

$$= \frac{1}{N_{k-1}} \sum_{j : \tilde{G}_{k-1}\left(V_{k-1}^j\right)=1} \left(1 - \tilde{M}_{(k-1):m}\left(V_{k-1}^j, \Delta\right)\right) \tag{9}$$

where (5) follows from the tower rule, and noting that $\tilde{p}_k$ is $\mathcal{G}_k$-measurable; (6) from the induction hypothesis; (7) from (4); (8) from expanding and integrating over $M_k(V_{k-1}^j, du)$; and (9) from the semigroup property.

Moreover, putting $k = 1$, we have

$$(9) = \frac{1}{N_{k-1}} \sum_{j=1}^{N_0} \left(1 - M_{0:m}\left(V_0^j, \Delta\right)\right).$$

Since $\left(1 - M_{0:m}\left(V_0^j, \Delta\right)\right) = \mathbb{P}\left(\tau_B < \tau_A | X(0)\right)$, we see that

$$\mathbb{E}\left[\tilde{p}\right] = \mathbb{E}\left[\prod_{i=1}^m \tilde{p}_m\right]$$
$$= \mathbb{P}(\tau_B < \tau_A)$$

as desired. □

**Remark 9** In this algorithm we have taken the simple choice to split paths of $\tilde{X}$ at the first (random) multiple of $T$ after which a crossing is guaranteed, $\tilde{\tau}_i$. Since this could result in a long gap between the crossing time and the splitting time, i.e. a large $\left(\tilde{\tau}_i - \tau_i\right)$, this may introduce some unwanted variance into the estimation and various other approaches to splitting could be implemented. One simple alternative is to choose in advance a reasonably fine deterministic time-grid, and to split at the first location on the time-grid after which crossing is guaranteed.

---

**Algorithm 6** Exact Multilevel Splitting

Given $\lambda$ together with $\tilde{G}_i, \tilde{M}_i$ for $i = 1, \ldots, m$, an initial number of particles $N_0$, and splitting ratios $R_1, \ldots, R_{m-1}$:

1. Initialise $S_1 = \cdots = S_m = \emptyset$.

2. For $j = 1, \ldots, N_0$:

   (a) Draw $\tilde{X}_0^j(0) \sim \lambda$, and simulate $\tilde{X}_0^j[0:T]$ together with its skeleton $\tilde{Z}_0^j$ as per Section 2.3.

   (b) Sample $(\tilde{G}_1(\tilde{Z}_1^j), \tilde{Z}_1^j) \sim$ Algorithm.4$(\tilde{Z}_0^j, A, B_1)$. If $\tilde{G}_1(\tilde{Z}_1^j) = 1$, add $\tilde{Z}_1^j$ to $S_1$.

3. Record $N_1 = |S_1|$.

4. For $i = 2, \ldots, m$:

   (a) if $N_{i-1} = 0$, return the estimate $\hat{p}^{\text{ex}} = 0$ of $p$.

   (b) Otherwise, given $S_{i-1} = \{\tilde{Z}_{i-1}^j\}_{j=1}^{N_{i-1}}$, for all pairs $\{(j,k)\}_{1 \le j \le N_{i-1}}^{1 \le k \le R_{i-1}}$ sample independently

   $$(\tilde{G}_i(\tilde{Z}_i^{j,k}), Z_i^{j,k}) \sim \text{Algorithm.4}(\tilde{Z}_{i-1}^j, A, B_i),$$

   and if $\tilde{G}_i(\tilde{Z}_i^{j,k}) = 1$, add $\tilde{Z}_i^{j,k}$ to $S_i$.

   (c) Set $N_i = |S_i|$.

5. Return

   $$\hat{p}^{\text{ex}} = \frac{N_m}{N_0 \prod_{i=1}^{m-1} R_i}.$$

---

The estimate $p^{\text{ex}}$ given by Exact Multilevel Splitting (Algorithm 6) is exactly the same as that given by Idealised Splitting with Coupling (Algorithm 5), since the particle systems defined by these algorithms are essentially identical.

**Corollary 1** $\mathbb{E}[\hat{p}^{\text{ex}}] = p$ for all $N \geq 1$.

Similarly to Algorithm 2, a Sequential Monte Carlo variant of Algorithm 6 may be constructed by replacing the splitting step with resampling: this is illustrated in Algorithm 7. Its advantages over Algorithm 6 are the same as the advantages of Algorithm 2 in Sect. 2.1, namely that there are no splitting ratios $R_i$ which need to be calibrated to ensure a stable particle system.

---

**Algorithm 7** Exact Multilevel Splitting via SMC

Given $\lambda$ together with $\tilde{G}_i, \tilde{M}_i$ for $i = 1, \ldots, m$, and a fixed number of particles $N$:

...

4(b)'. Otherwise, given $S_{i-1} = \{\tilde{Z}_{i-1}^j\}_{j=1}^{N_{i-1}}$, sample $a_1, \ldots, a_N$ independently and uniformly at random (with replacement) from $1, \ldots, N_{i-1}$, and sample independently $(\tilde{Z}_i^j, \tilde{G}_i(\tilde{Z}_i^j)) \sim$ Algorithm.4$(\tilde{Z}_{i-1}^{a_j}, A, B_i)$.

...

5. Estimate

$$\hat{p}^{\text{SMC}} = \frac{\prod_{i=1}^m N_i}{(N)^m} = \prod_{i=1}^m \left(\frac{N_i}{N}\right).$$

---

# 4 Illustrations

The examples in this section were carried out using a single core on an Intel Xeon E5-2440 processor with an advertised clock speed of 2.40GHz.

## 4.1 Brownian Motion in One Dimension

Our first illustrative example uses the $\varepsilon$-strong scheme for Brownian motion of Beskos et al. (2012); Pollock et al. (2016), in a setting in which the exact solution is known. In one dimension, the reaction co-ordinate may be taken to be the identity function. We choose $A = (-\infty, 0]$, $B = [3^{18}, \infty)$, $B_i = [3^i, \infty)$ for $i = 1, \ldots, 17$, with initial point $x_0 = 1$. It is well-known that for real $0 < a < b$, the probability that a Brownian path started at $a$ reaches $b$ before 0 is $a/b$ — as can be verified with a simple optional stopping argument. Therefore the target probability is $3^{-18} \approx 2.58 \times 10^{-9}$.

The $\varepsilon$-strong scheme in question has some additional features which allow a substantial improvement in speed to that given in Algorithm 7 in this exceptionally simple setting. Over a given time interval $[s, t]$, for a Brownian motion $(X(r))_{r \in [s,t]}$ it is possible to sample bounds

$$L^{\downarrow}(s,t) \leq \inf_{r \in [s,t]} X(r) \leq L^{\uparrow}(s,t) \, , \, U^{\downarrow}(s,t) \leq \sup_{r \in [s,t]} X(r) \leq U^{\uparrow}(s,t).$$

An $\varepsilon$-strong algorithm over $[0, T]$ in the sense given in Sect. 2.3 can be recovered by choosing partitioning $[0, T]$ into suitably small time intervals $[s_i, t_i]$, and taking

$$\tilde{X}^{\varepsilon}(r) = \frac{1}{2} \left( U^{\uparrow}(s_i, t_i) - L^{\downarrow}(s_i, t_i) \right)$$

for $r \in [s_i, t_i]$. This corresponds to taking $C_i = [L^{\downarrow}(s_i, t_i), U^{\uparrow}(s_i, t_i)]$ in the notation of Sect. 3.1. This does not make use of the extra information available in $(L^{\uparrow}, U^{\downarrow})$. In particular, to assess whether $X$ has crossed above the point $b$ over $[s_i, t_i]$, it is sufficient to find that $U^{\downarrow}(s_i, t_i) > b$, since this guarantees the maximum of $X$ is large enough. This is easier to check than the stricter condition that $C_i \subset (b, \infty)$, or equivalently that $L^{\downarrow}(s_i, t_i) > b$.

We compare the exact estimator to Euler-Maruyama-type schemes (see (2)) with three levels of resolution. The initial step sizes for the schemes are taken to be 0.01, 0.005 and 0.001. In this example we exploit the simplicity of the problem at hand and the time-scaling property of Brownian motion to allow the Euler-Maruyama scheme to maintain a constant level of relative error over time: in each scheme, when level $B_i$ is reached the step size is multiplied by $3^2 = 9$. This scaling which depends upon analytical techniques which would not be available in more realistic problems was essential in order to achieve a reasonable calculation time (this choice ensures the expected number of Euler-Maruyama steps until a crossing is decided remains constant as the scale of the problem grows). 500 estimators were produced for each procedure, each with population of 1000 particles. The results are shown in Fig. 5.

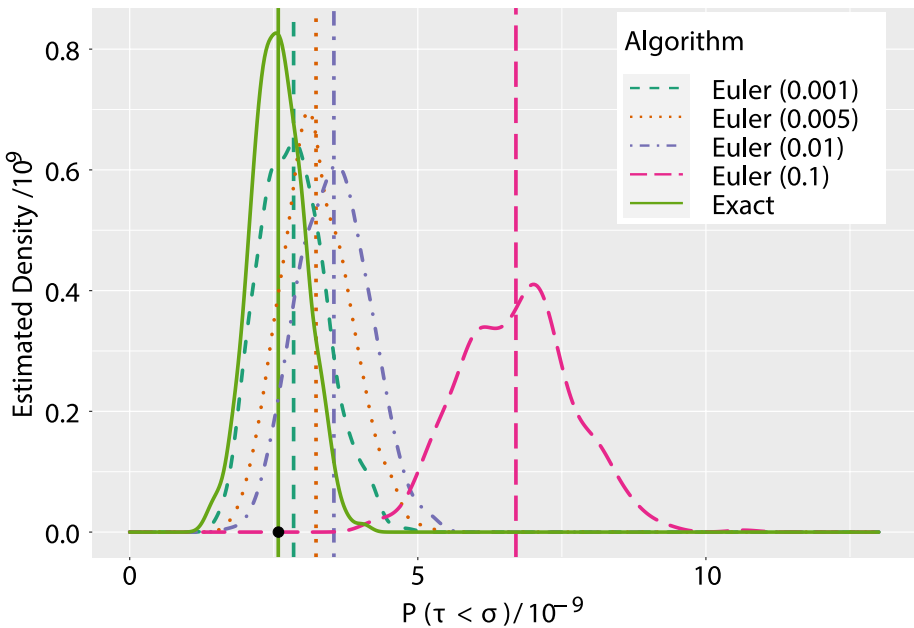

**Fig. 5** Kernel density estimates for the various schemes. The mean of each density estimate is shown as a coloured vertical line. The black dot on the ordinal axis indicates the true probability

Typical run-times for a single sample were 19s for the exact-MLS algorithm, 32s for the Euler-Maruyama-0.005 scheme, 157s for the Euler-Maruyama-0.001 scheme. This demonstrates that in favourable circumstances exact MLS can yield estimates of rare event probabilities at significantly lower computational cost than that at which discretisation-based methods can reach an acceptable level of bias. In other settings the cost of exact methods can be somewhat higher, as the next example will demonstrate.

### 4.2 A Bivariate Example

Our second example illustrates Algorithm 7 in a pure form, is a two-dimensional problem. The random process is again taken to be Brownian motion initialised at $W_0 = \left(\frac{1}{2}, \frac{1}{2}\right)$. The reaction co-ordinate is chosen to be $\xi(x, y) = \min(x, y)$, and the levels are chosen to be $A = \xi^{-1}((-\infty, 0)), B = \xi^{-1}\left(\left(2^{\frac{20}{2}}, \infty\right)\right), B_i = \xi^{-1}\left(2^{(i/2+1)}, \infty\right)$ for $i = 1, \ldots, 18$. We are not aware of any simple means by which the rare event probability can be analytically obtained in this case.

The $\varepsilon$-strong algorithm used is that of Beskos et al. (2012); Pollock et al. (2016), in the same fashion as Sect. 4.1. Again, the exact estimator is compared to three Euler-Maruyama schemes of (2) with increasing degrees of fineness, in this case with initial step-sizes 0.1, 0.05 and 0.01. The step-sizes were again re-scaled at each new level to ensure an approximately constant relative error, this time by a factor of $2^2 = 4$. For each scheme, 250 trials were simulated, each using 100 particles. Typical running times were 13s for the Euler-Maruyama-0.1 scheme, 21s for the Euler-Maruyama-0.05 scheme, 1m 45s for the Euler-Maruyama-0.05 scheme and 534m for the exact scheme.

Although the running time for exact MLS is significantly longer than that of the discrete schemes in this case, we believe that a more careful attempt to tune and adapt its parameters could substantially reduce the difference. However, the optimal choice of parameters will depend on the particular application, and since our aim is to provide proof-of-concept validation of the generic methodology developed in this paper we have not attempted to do so.

Figure 6 shows a kernel density estimate of the sampling distribution of the resulting estimators obtained using the `geom_density` function of ggplot2 (Wickham, 2016), using its default choice of bandwidth. As the Euler-Maruyama scheme increases in fineness, the resulting estimate appears closer to the estimate from exact MLS.

## 5 Discussion

We have presented the first algorithm for the exact estimation of a class of rare event probabilities for continuous-time processes. Our method has been to directly replace discrete approximations of continuous-time sample paths with $\varepsilon$-strong samples of the same paths, in order to obtain the sequence $\hat{p}_i$ of conditionally unbiased estimates required for multilevel splitting.

There is considerable ongoing effort in the development of $\varepsilon$-strong simulation methods for a broad class of stochastic processes, and their application (see for instance Mider et al. (2019)). It is likely that further development in this area will allow the approach described within this paper to be applied with greater efficiency to a broader class of stochastic processes. Naturally, it is likely that the exactness of this approach will always come at the expense of a computational cost that exceeds that of discretisation-based schemes. This cost can be partially offset against the need to assess the bias inherent in such schemes.

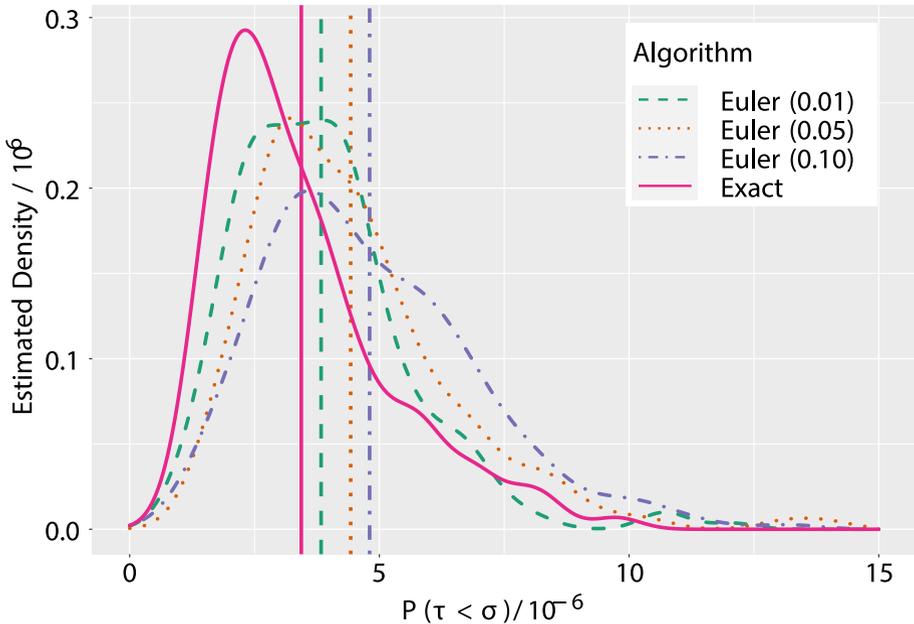

**Fig. 6** Kernel density estimates for exact MLS and for an Euler-Maruyama-discretisation-based method with three different discretisation step sizes. The mean of each density estimate is shown as a coloured vertical line

One important assumption that we have made is the ability to calculate infima and suprema of the reaction co-ordinate over the sets $C_k$ which constrain a sample path $X$. The problem of assessing whether these $C_k$ intersect the MLS sets $B_i$ has links to the problem of *collision detection* studied in computer graphics (see eg Kockara et al. (2007)). Insight from this field might allow a more careful classification of suitable $\xi$ for a given $C_k$, or for the assumption to be weakened in certain circumstances.

The *alive particle filter* of Del Moral et al. (2015) provides an alternative approach to implementing SMC algorithms with $\{0, 1\}$-valued potential functions which also provide unbiased estimates of normalizing constants and thus, in the present context, of the rare event probability. We did not experience difficulties with extinction of the particle system in our experiments, but incorporating that approach within the exact MLS framework that we present would be interesting because it would automatically mitigate the influence of poorly chosen intermediate levels, albeit at the cost of further randomizing the computational cost.

As we noted in Sect. 2.1, one strength of this method is that other quantities relating to the distribution of *paths* which reach the rare event set $B$ may also be estimated. These require the sort of detailed information given by $\varepsilon$-strong simulation. However, if the rare event probability is the only quantity of interest, it is possible that approaches to obtaining unbiased estimates of rare event probabilities *without* requiring the full machinery of unbiased MLS implementations are sufficient.

The construction of computable unbiased estimators via a sequence of asymptotically biased estimators has received much attention in recent years, following Glynn and Rhee (2014). This technique has been directly extended to estimating expectations of functionals of SDE paths in Rhee and Glynn (2015), and recently this approach has been further extended to non-linear filtering problems in Jasra et al. (2020). In ongoing work we have started to explore the possibility that a related approach could allow for exact rare event estimation in a more general setting, and with greater scope for practical applications.

**Funding** This work was supported by the Engineering and Physical Sciences Research Council under grant numbers EP/L016710/1, EP/R034710/1 and EP/T004134/1 and by two Alan Turing Institute programmes; the Lloyd's Register Foundation programme on 'Data-centric engineering' and the UK Government's 'Defence and security' programme.